\documentclass[12pt,twoside]{amsart}

\usepackage{amssymb,latexsym,bbm,graphicx,epsfig,epic,eepic,oldgerm,psfrag}
\usepackage{a4wide}

\usepackage{xypic} 
\input xy

\xyoption{all}

\theoremstyle{plain}
\newtheorem{theorem}{Theorem}[section]
\newtheorem{lemma}[theorem]{Lemma}
\newtheorem{proposition}[theorem]{Proposition}
\newtheorem{corollary}[theorem]{Corollary}
\newtheorem*{remark*}{Remark}
\newtheorem*{remarks*}{Remarks}
\newtheorem{remark}[theorem]{Remark}
\newtheorem{remarks}[theorem]{Remarks}

\newtheorem*{example*}{Example}
\newtheorem*{examples*}{Examples}
\newtheorem{definition}[theorem]{Definition}

\newcommand{\proofend}{\hspace*{\fill} $\Box$\\}

\def\1{\:\!}
\def\2{\;\!}
\def\s{\smallskip}
\def\m{\medskip}

\def\Volg{\operatorname {Vol}}

\def\Diffc0{\operatorname{Diff^c_0}}

\def\Sympc0{\operatorname{Symp^c_0}}

\def\rank{\operatorname{rank}}

\def\top{\operatorname{top}}

\def\vol{\operatorname{vol}}
\def\Aut{\operatorname{Aut}}

\def\SL2{\operatorname{SL_2}}

\def\ga{\alpha}

\def\gg{\gamma}

\def\gf{\varphi}

\def\gs{\sigma}

\def\cc{{\mathcal C}}

\def\cp{{\mathcal P}}

\def\cs{{\mathcal S}}

\def\FF{\mathbbm{F}}

\def\NN{\mathbbm{N}}

\def\QQ{\mathbbm{Q}}
\def\RR{\mathbbm{R}}

\def\ZZ{\mathbbm{Z}}

\def\fm{{\mathfrak m}}

\def\pp{\partial}

\def\pr{{\rm pr}}

\def\ra{\rightarrow}

\def\ni{\noindent}

\def\m{\medskip}

\def\proof{\noindent {\it Proof. \;}}

\begin{document}

\title[]{Healthy vector spaces and spicy Hopf algebras \\ 
(with applications to the growth rate of  geodesic 
chords and to intermediate volume growth 
on manifolds of non-finite type)}

\author{Urs Frauenfelder}
\thanks{UF partially supported by the Basic Research fund 2013004879 of the
Korean government}
\address{
    Urs Frauenfelder,
    Department of Mathematics and Research Institute of Mathematics,
    Seoul National University}
\email{frauenf@snu.ac.kr}

\author{Felix Schlenk}  
\thanks{FS partially supported by SNF grant 200020-144432/1.}
\address{Felix Schlenk,
Institut de Math\'ematiques,
Universit\'e de Neuch\^atel}
\email{schlenk@unine.ch}

\date{\today}

%

\begin{abstract}
We give lower bounds for the growth of the number of Reeb chords and for
the volume growth
of Reeb flows on spherizations over closed manifolds~$M$
that are not of finite type,
have virtually polycyclic fundamental group, 
and satisfy a mild assumption on the homology of the based loop space.
For the special case of geodesic flows, these lower bounds are:

\s
\begin{itemize}
\item[(i)]
For any Riemannian metric on~$M$, any pair of non-conjugate points $p,q \in M$,
and every component~$\cc$ of the space of paths from~$p$ to~$q$,
the number of geodesics in~$\cc$ of length at most~$T$
grows at least like $e^{\sqrt T}$.

\s  
\item[(ii)]                
The exponent of the volume growth of any geodesic flow on~$M$ is at least $1/2$. 
\end{itemize}

\m \ni
We obtain these results by combining new algebraic results 
on the growth of certain filtered Hopf algebras
with known results on Floer homology.
\end{abstract}

\maketitle

\section{Introduction and main results}

\ni
Consider a closed connected manifold~$M$. 
Then $M$ is said to be of {\it finite type}\/ if its universal cover~$\widetilde M$
is homotopy equivalent to a finite CW-complex.
Equivalently, there is $k \in \{2, \dots, \dim M \}$ such that $H_k(\widetilde M)$
is not finitely generated. Let $\fm = \fm (M)$ be the minimal such~$k$.


A finitely generated group $G$ is called {\it polycyclic}\/ if it
admits a subnormal series with cyclic factors.
Moreover, $G$ is {\it virtually polycyclic}\/ if it has a polycyclic subgroup 
of finite index.

Choose a point $p \in M$ and let $\Omega_0M$ be the space of contractible continuous loops in~$M$ based at~$p$. 
For $T>0$ denote by $\Omega_0^T M$ the space of contractible piecewise smooth loops based at~$p$ whose length (with respect to a fixed Riemannian metric on~$M$) 
is at most~$T$. 
The inclusions $\iota^T \colon \Omega_0^T M \hookrightarrow \Omega_0 M$ induce the maps 
$\iota^T_* \colon H_*(\Omega_0^T M) \to H_*(\Omega_0 M)$ in homology.
Note that the image $\iota_*^T H_*(\Omega_0^T M)$ is the part of the homology of~$\Omega_0M$
that is generated by cycles made of piecewise smooth loops of length~$\le T$.

We say that a real-valued function~$f$ defined on~$\NN$ or on~$\RR_{>0}$ 
{\it grows at least like $e^{\sqrt T}$}
if there exists a constant~$c>0$ such that 
$f(T) \ge c \2 e^{\sqrt T}$ for all large enough~$T$.

\begin{theorem} \label{t:main}
Let $M$ be a closed connected manifold that is not of finite type
and has virtually polycyclic fundamental group.
Assume also that there exists a field~$\FF$ such that $H_{\fm}(\widetilde M;\FF)$ is 
infinite-dimensional.
Then the function of~$T$
$$
\dim \iota_*^T H_* (\Omega_0^T M; \FF) 
$$
grows at least like $e^{\sqrt T}$.
\end{theorem}

\begin{remarks} \label{rem:HP}
{\rm
{\bf 1.}
The assumption that $H_{\fm} (\widetilde M;\FF)$ is infinite-dimensional
is equivalent to the assumption that $\pi_{\fm}(M) \otimes \FF$ is infinite-dimensional. 
 
\s
{\bf 2.}
We conjecture that the rank of $\iota_*^T H_* (\Omega_0^T M)$ grows exponentially for every
manifold~$M$ of non-finite type 
(cf.\ the Question in~\cite[p.\ 289]{PatPet06} and the discussion in~\cite[\S 7.1]{FraLabSch13}).
Theorem~\ref{t:main} proves ``not quite half'' of this conjecture.

\s
{\bf 3.}
Along the proof we shall see that 
{\it $H_{\fm-1}(\Omega_0M;\FF)$ is infinite-dimensional and that
$$
\dim H_{(\fm-1)n} (\Omega_0M;\FF) \,\ge\, q(n) \quad \mbox{for }\, n \ge 2,
$$
where $q(n)$ is the number of partitions of $n$ into distinct parts.}
The fact that this function grows like $e^{\sqrt n}$ explains the lower bound
in~Theorem~\ref{t:main}.

\s
{\bf 4.}
Our assumption that there exists a field~$\FF$ such that $H_\fm (\widetilde M;\FF)$
is infinite-dimensional plays an important role in our proof. 
This is illustrated by an example in Section~\ref{ss:exa}. 
In Section~\ref{ss:exa.neu} we give a class of examples that meet the assumptions of Theorem~\ref{t:main}.
}
\end{remarks}

\subsection{Applications}
Lower bounds for the rank of the homology of the sublevel sets $\Omega^T M$ are of interest because they classically lead, by Morse theory, to lower bounds for the number of geodesics of length $\le T$ between non-conjugate points.
Somewhat less classically, they also lead to lower bounds for the topological entropy 
of geodesic flows. 
Moreover, adding the tool of Floer homology, one gets lower bounds for the number 
of Reeb chords and for the topological entropy of Reeb flows on spherizations.

Before stating our two corollaries, we briefly recall 
what Reeb flows on spherizations are. 
Details can be found in the introduction of~\cite{FraLabSch13}.

\m \ni
{\bf Reeb flows on spherizations.}
Consider a closed manifold~$M$.
The positive real numbers~$\RR_+$ freely act on the cotangent bundle~$T^*M$
by $r\1 (q,p) = (q,r\1 p)$. While the canonical $1$-form $\lambda = p\1 dq$ on $T^*M$ does not 
descend to the quotient $S^*M := T^*M / \RR_+$, 
its kernel does and defines a contact structure~$\xi$ on~$S^*M$.
We call the contact manifold $(S^*M,\xi)$ the {\it spherization}\/ of~$M$.
This contact manifold is co-orientable. 
The choice of a nowhere vanishing 1-form~$\alpha$ on~$S^*M$ with $\ker \alpha = \xi$
(called a {\it contact form}) defines a vector field $R_\alpha$ 
(the {\it Reeb vector field }\/of~$\alpha$) 
by the two conditions
$d\alpha (R_\alpha, \cdot ) = 0$, $\alpha (R_{\alpha}) = 1$.
Its flow $\varphi_\alpha^t$ is called the {\it Reeb flow}\/ of~$\alpha$. 

To give a more concrete description of the manifold $(S^*M,\xi)$ 
and the flows $\gf_\alpha^t$,
consider a smooth hypersurface~$\Sigma$ in $T^*M$ which is {\it fiberwise starshaped}
with respect to the zero-section: 
For every $q \in M$ the set $\Sigma_q := \Sigma \cap T_q^*M$ bounds a set in~$T_q^*M$ 
that is strictly starshaped with respect to the origin of~$T_q^*M$.
The hyperplane field $\xi_\Sigma := \ker (\lambda |_\Sigma)$ is a contact structure on~$\Sigma$,
and the contact manifolds $(S^*M,\xi)$ and $(\Sigma, \xi_\Sigma)$ are isomorphic.

Let $\gf_\Sigma^t$ be the Reeb flow on $\Sigma$ defined by the contact form $\lambda |_\Sigma$.
The set of Reeb flows on $(S^*M,\xi)$ can be identified with the Reeb flows $\gf_\Sigma^t$ 
on the set of fiberwise starshaped hypersurfaces~$\Sigma$ in $T^*M$.
The flows $\gf_\Sigma^t$ are restrictions of Hamiltonian flows:
Consider a Hamiltonian function $H \colon T^*M \to \RR$ such that 
$\Sigma = H^{-1}(1)$ is a regular energy surface and such that $H$ is fiberwise homogeneous of degree~one
near $\Sigma$.
For the Hamiltonian flow $\gf_H^t$ we then have $\gf_H^t |_\Sigma = \gf_\Sigma^t$.
It follows that geodesic flows and Finsler flows 
(up to the time change $t \mapsto 2t$)
are examples of Reeb flows on spherizations.
Indeed, for geodesic flows the $\Sigma_q$ are ellipsoids, 
and for (symmetric) Finsler flows the $\Sigma_q$ are (symmetric and) convex. 
The flows $\gf_\Sigma^t$ for varying~$\Sigma$ are very different, in general, 
as is already clear from looking at geodesic flows on a sphere.
In this paper we give uniform lower bounds for the growth of Reeb chords and for the complexity 
of all these flows on~$(S^*M,\xi)$ for manifolds~$M$ as in Theorem~\ref{t:main}.

\m \ni
{\bf Growth of Reeb chords between two fibers.}
We say that $p,q \in M$ are {\it non-conjugate points}\/ of the Reeb flow~$\gf_\alpha$
if $\bigcup_{t>0} \gf_\ga^t (S^*_pM)$ is transverse to~$S^*_qM$.
Given $p \in M$, the set of $q \in M$ that are non-conjugate to~$p$ has full measure in~$M$ 
by Sard's theorem. This notion of being non-conjugate generalizes the one in  
Riemannian geometry (defined in terms of Jacobi fields).

For $p,q \in M$ denote by $\cp_{pq}$ the space of continuous paths in~$S^*M$ from~$S^*_pM$ to~$S^*_qM$, 
and by $\Omega_{pq}M$ the space of continuous paths in~$M$ from~$p$ to~$q$.
We shall assume throughout that $\dim M \ge 3$, since otherwise $M$ is of finite type.
The fibers $S^*_qM$ of the projection $pr \colon S^*M \to M$ are then simply connected, and so 
$pr$ induces a bijection on the components of $\cp_{pq}$ and $\Omega_{pq}M$.
The space $\Omega_{pq}M$ is homotopy equivalent to $\Omega_pM := \Omega_{p,p}M$, whose components
are parametrized by the elements of the fundamental group $\pi_1(M,p)$.

\begin{corollary} \label{c:chords}
Assume that $M$ is not of finite type, has virtually polycyclic fundamental group,
and that there is a field~$\FF$ such that $H_{\fm}(\widetilde M;\FF)$ is infinite-dimensional.
Let $(S^*M, \xi)$ be the spherizations of~$M$.
Then for any Reeb flow~$\gf_\ga$ on~$(S^*M, \xi)$, any pair of non-conjugate points $p,q \in M$
and every component~$\cc$ of~$\Omega_{pq}M$,
the number of Reeb chords from~$S_p^*M$ to~$S_q^*M$ that belong to~$\cc$ grows in time
at least like $e^{\sqrt T}$.
\end{corollary}

\begin{remarks}
{\rm
{\bf 1.}
For the special case of geodesic flows, the time parameter equals the length 
run through.
Hence the corollary translates to assertion~(i) of the abstract.

To illustrate the corollary,
we choose a Riemannian metric on~$M$ and a point $p \in M$.
Let $C(p)$ be the cut locus of~$p$. The subset $M \setminus C(p)$ is diffeomorphic
to an open ball~\cite{Ozo76} and has full measure in~$M$.
For every $q \in M \setminus C(p)$ there is a unique shortest geodesic~$c_q$ from~$p$ to~$q$.
Call a path $\gg \in \Omega_{pq}M$ {\it contractible}\/ if $c_q^{-1} \circ \gg$ 
is contractible in~$\Omega_pM$.
The set~$U_p$ of points in $M \setminus C(p)$ that are not conjugate to~$p$ is also of full 
measure in~$M$. Under the hypothesis of Corollary~\ref{c:chords},
for every $q \in U_p$ the number of contractible geodesics from~$p$ to~$q$ 
of length~$\le T$ grows at least like $e^{\sqrt T}$.

\s
{\bf 2.} 
Virtually polycyclic groups are either virtually nilpotent or have exponential
growth~\cite{Wol68}.
If the fundamental group~$\pi_1(M)$ of a closed manifold~$M$ has exponential growth, 
then the number of Reeb chords from $S_p^*M$ to~$S_q^*M$ grows exponentially in time
for any, possibly conjugate, pair of points~$p,q$
(see \cite[Corollary~1]{MacSch11} for Reeb flows).
Indeed, one finds one Reeb chord for each element of~$\pi_1M$.
In Corollary~\ref{c:chords}, however, we find ``$e^{\sqrt T}$ many'' Reeb chords for each element of~$\pi_1M$.

If one is only interested in the growth of Reeb chords from $S_p^*M$ to~$S_q^*M$,
without specifying the component~$\cc$, then Corollary~\ref{c:chords}
is interesting only for virtually nilpotent fundamental groups,
which by Gromov's theorem from~\cite{Gro81} are exactly those fundamental groups
that grow polynomially.
Indeed, it is believed that every finitely presented group that grows more than 
polynomially grows exponentially~\cite[Conjecture~11.3]{Gri08},
and even for finitely generated groups of intermediate growth it is believed
that they must grow at least like $e^{\sqrt T}$, cf.\ \cite{Gri12}. 

\s
{\bf 3.}
Let $\mu(\gg)$ be the Conley--Zehnder index of a non-degenerate Reeb chord~$\gg$
on~$(S^*M;\ga)$,
normalized such that for geodesic flows $\mu(\gg)$ is the Morse index 
of the non-degenerate geodesic~$\gg$
(i.e.\ the number of conjugate points, counted with multiplicities, along~$\gg$).
In the situation of Corollary~\ref{c:chords}, Remark~\ref{rem:HP}.3 shows that
for every component of~$\Omega_{pq}$
the number of Reeb chords from~$S^*_pM$ to~$S^*_qM$ of index $\mu (\gg) = k$ 
is infinite for $k=\fm-1$ and at least $q(n)$ if $k=(\fm-1)n$ and $n \ge 2$.
}
\end{remarks}

\ni
{\bf Intermediate volume growth.}
%
%
Consider a smooth diffeomorphism $\gf$ of a closed manifold~$X$.
Denote by $\cs$ the set of smooth compact submanifolds of~$X$.
Fix a Riemannian metric~$\rho$ on~$X$,
and denote by $\Volg (\gs)$ the $j$-dimensional volume 
of a $j$-dimensional submanifold $\gs \in \cs$
computed with respect to the measure on~$\gs$ induced by~$\rho$. 
For $a \in (0,1]$
define the {\it intermediate volume growth}\, of $\gs \in \cs$ by 
\begin{equation} \label{def:slowvol}
\vol^a (\gs; \gf) \,=\, 
 \liminf_{n \ra \infty} 
 \frac{\log \Volg \left( \gf^n ( \gs ) \right)}{n^a}  \,\in\, [0,\infty] ,
\end{equation}
and define the {\it intermediate volume growth}\, of $\gf$ by
$$                  
\vol^a (\gf) \,=\, \sup_{\gs \in \cs} \vol^a (\gs;\gf) \,\in\, [0,\infty] .
$$
Notice that these invariants do not depend on the choice of~$\rho$.
Finally define the {\it volume growth exponent}\, of~$\gf$ by
$$
\exp_{\vol} (\gf) \,:=\, 
\inf \left\{ a \mid \vol^a(\gf) < \infty \right\} .
$$
Thus $\exp_{\vol} (\gf)$ is ``the largest $a \in [0,1]$ such that some submanifold grows under~$\gf$ like~$e^{n^a}$.''
Since $\vol^1(\gf) \le (\dim X) \max_{x  \in X} \| D \gf (x) \| < \infty$ we have 
$\exp_{\vol} (\gf) \in [0,1]$.
The intermediate volume growth and the volume growth exponent of a smooth flow~$\gf^t$ on~$X$ are defined as $\vol^a (\gf^1)$ and $\exp_{\vol} (\gf^1)$.

\begin{remark}
{\rm
By a celebrated result of Yomdin~\cite{Yom87} and Newhouse~\cite{New88}, 
the volume growth $\vol^1 (\gf)$ agrees with the topological entropy $h_{\top} (\gf)$.
Proceeding as above define for $a \in (0,1]$ 
the {\it intermediate topological entropy} $h_{\top}^a (\gf)$.
We unfortunately do not know whether $\vol^a (\gf) = h_{\top}^a(\gf)$ 
or at least $\vol^a (\gf) \le h_{\top}^a(\gf)$ 
also for $a \in (0,1)$.
}
\end{remark}

Uniform lower bounds for the volume growth or the topological entropy of geodesic flows
were found in~\cite{Din71, Pat97, Pat:book, PatPet06}, 
and these results were generalized to Reeb flows in~\cite{MacSch11} and, on a polynomial scale,
in~\cite{FraLabSch13}.
Results in~\cite{FraLabSch13, MacSch11} show that {\it the volume growth exponent
$\exp_{\vol}(\gf_\ga)$ is bounded from below by the maximum of the growth exponent 
of the function $n \mapsto \rank \iota_*^T H_* (\Omega_0^T M)$ and the growth exponent 
of the growth function of the fundamental group $\pi_1(M)$.}
In particular, for manifolds with fundamental group of exponential growth it is shown in~\cite{MacSch11} 
that $\vol^1(\gf_\ga) >0$ for any Reeb flow $\gf_\ga$ on~$S^*M$.
Since virtually polycyclic groups are either virtually nilpotent or have exponential 
growth, we therefore restrict ourselves now to manifolds with virtually nilpotent 
fundamental group.
Since all other fundamental groups of closed manifolds are believed to have 
exponential growth, this is a minor hypothesis.

\begin{corollary} \label{c:ent}
Assume that $M$ is not of finite type, has virtually nilpotent  
fundamental group,
and that there is a field~$\FF$ such that $H_{\fm}(\widetilde M;\FF)$ is infinite-dimensional.
Let $(S^*M, \xi)$ be the spherizations of~$M$.
Then 
$$
\vol^{1/2} (S_p^*M;\gf_\ga) >0
$$
for every fiber $S_p^*M$ of~$S^*M$ and every Reeb flow~$\gf_\ga$ on~$(S^*M, \xi)$.
In particular, 
$\vol^{1/2} (\gf_\ga) >0$ and $\exp_{\vol}(\gf_\ga) \ge 1/2$ for every Reeb flow~$\gf_\ga$ on~$(S^*M, \xi)$.
\end{corollary}

\ni
{\bf The method.}
We end this introduction with comparing our approach to previous approaches.
As mentioned earlier, lower bounds for the homology of the sublevel sets~$\Omega^T M$ 
follow easily from the growth of the fundamental group~$\pi_1(M)$, 
since its growth is the growth of~$H_0(\Omega^TM)$.
If $\pi_1(M)$ is finite, Gromov found a way to bound the rank of~$H_*(\Omega^TM)$ 
from below by the rank of $\oplus_{i=0}^{cT} H_i(\Omega M)$, where $c$ is a constant depending only 
on the Riemannian metric.
His ingenious argument is purely geometric, see~\cite{Gro78, Gro07, Pat97}, 
and uses that~$\widetilde M$ is compact.
In~\cite{PatPet04,PatPet06} Paternain--Petean generalized Gromov's construction 
to manifolds with infinite fundamental group, 
by mapping simply connected complexes~$K$ with rich loop space homology into~$\widetilde M$.
This method gives good lower bounds for the rank of~$H_*(\Omega^T M)$ if one can find 
such complexes and mappings $f \colon K \to \widetilde M$ for which the homology of 
$\Omega f (\Omega K)$ still has large rank. 
This works well for many manifolds, e.g.\ for most connected sums, 
for manifolds of finite type with $\pi_*(\Omega M) \otimes \QQ$ infinite-dimensional,
and in small dimensions, see~\cite{PatPet04, PatPet06}.
For general manifolds of non-finite type, however, we were not even able to prove
linear growth of the rank of $H_*(\Omega^T M)$ by this method.
Citing \mbox{G.\ Paternain},  
``it is as if one has so much topology that it becomes unmanageable.''

Our approach to finding lower bounds for the rank of~$H_*(\Omega^TM)$ for manifolds of non-finite type
is neither geometric nor topological, but algebraic. 
The main point is that we use the Hopf algebra structure of $H_*(\Omega \widetilde M;\FF)$ 
over the group ring~$\FF [\pi_1(M)]$ for a suitable field~$\FF$, 
and prove a Poincar\'e--Birkhoff--Witt type theorem for this algebra.
These algebraic results (specifically Theorems~\ref{t:healthy} and~\ref{t:spicy} 
and Proposition~\ref{p:quantpbw}) are the main findings of this paper. 
They cover Sections~\ref{s:healthy} to~\ref{s:PBW}.
Theorem~\ref{t:main}, that is proven in Section~\ref{s:proof}, 
follows readily from these algebraic results, 
and Corollaries~\ref{c:chords} and~\ref{c:ent}, that are proven in Section~\ref{s:cor}, 
follow from Theorem~\ref{t:main} in the usual way.

\m
\ni
{\bf Acknowledgments.}
We wish to thank Kenji Fukaya for suggesting to us to bring into play
Serre's Hurewicz theorem.
The present work is part of the author's activities within CAST, 
a Research Network Program of the European Science Foundation.


\section{Healthy $G$-vector spaces} \label{s:healthy}

\ni
Let $G$ be a group, $V$ a vector space over the field~$\FF$, and
$\rho \colon G \to \Aut (V)$ a representation of~$G$. 
Since the representation is fixed throughout our discussion, 
we abbreviate $gv=\rho(g)v$ for $g \in G$ and $v \in V$, 
and we consider the $G$-vector space~$V$ as a module over the group ring~$\FF G$.
A group $G$ is called \emph{polycyclic}\/ if there exists a subnormal series
$$
1=G_0 \lhd G_1 \lhd \cdots \lhd G_n=G
$$
such that all factors $G_i/G_{i-1}$, $1 \leq i \leq n$ are cyclic.
A group $G$ is called a \emph{virtually polycyclic}\/ if $G$ has a polycyclic subgroup 
of finite index. 
Eleven characterisations of polycyclic groups are given in~\cite[Proposition~4.1]{Wol68}.
In this section we prove the following result.  

\begin{theorem} \label{t:healthy}
Assume that $G$ is a virtually polycyclic group acting linearly 
on an infinite-dimensional vector space~$V$ 
such that $V$ viewed as an $\FF G$-module is finitely generated. 
Then there exists $v \in V$ and $g \in G$ such that the 
sequence of vectors $\left( g^i v \right)_{i \in \ZZ}$ is linearly independent.
\end{theorem}

Before embarking on the proof of the theorem we introduce some notation. 
For the following discussion it is irrelevant that $G$ is virtually polycyclic. 
Given $g \in G$ and $v \in V$ we denote the subvector space of~$V$ spanned by the vectors $g^iv$ by
$$
W_v^g \,:=\, \big\langle g^i v \mid i \in \ZZ \big\rangle \,\subset\, V .
$$

\begin{definition}
{\rm
A vector $v \in V$ is \emph{healthy}\/ if there exists $g \in G$ such that $W_v^g$ 
is infinite-dimensional. 
A vector $v \in V$ which is not healthy is called \emph{sick}.
}
\end{definition}

The reason why we are interested in finding healthy vectors is the following observation from 
linear algebra.

\begin{lemma}\label{linalg}
If $g \in G$ and $v \in V$, then the following are equivalent. 

\s
\begin{itemize}
\item[(i)] $W^g_v$ is infinite-dimensional.

\s
\item[(ii)] The sequence of vectors $\left( g^iv \right)_{i \in \ZZ}$ is linearly independent.
\end{itemize}
\end{lemma}

\proof
That (ii) implies (i) is clear. 
It remains to show that if (ii) does not hold, then $W^g_v$ is finite-dimensional. 
The case $v=0$ is trivial as well, so we assume that $v \neq 0$. 
Since then $g^{-i}(g^i v) = v \neq 0$, we conclude that $g^i v \neq 0$ for every $i \in \ZZ$. 
Hence, if (ii) does not hold, there exist $m<n \in \ZZ$ and scalars $a_i \in \FF$ 
for $m \leq i \leq n-1$ such that
\begin{equation}\label{lindep}
g^nv=\sum_{i=m}^{n-1}a_i \2 g^i v, \quad a_m \neq 0.
\end{equation}
By applying $g^{-m}$ to \eqref{lindep} we can assume without loss of generality that~$m=0$.
Applying $g$ and $g^{-1}$ to~\eqref{lindep} we obtain inductively that
\begin{equation}\label{gensys}
W^g_v = \big\langle g^i v \mid 0 \leq i \leq n-1 \big\rangle.
\end{equation}
This shows that $W^g_v$ is finite-dimensional. The lemma follows. 
\proofend

If $v \neq 0$ is a sick vector, we can define in view of~\eqref{gensys} the function
$d_v \colon G \to \NN$ by
\begin{equation}\label{dv}
d_v(g) := \min \Bigl\{ n \in \NN \mid W^g_v = \langle g^i v \mid 0 \leq i \leq n-1 \rangle \Bigr\}
\,=\, \dim W^g_v .
\end{equation}

\begin{definition}
{\rm
The $G$-vector space $V$ is called \emph{sick}\/ if all its vectors are sick. 
A $G$-vector space $V$ which is not sick is called \emph{healthy}, i.e., $V$ contains a healthy vector.
}
\end{definition}

Note that a healthy vector space still contains sick vectors. 
Indeed, the zero vector is always sick. 
Moreover, observe that the concept of a healthy $G$-vector space is only of interest 
if both the cardinality of the group and the dimension of the vector space are infinite, 
since otherwise $V$ is automatically sick. 

Although we fix the representation~$\rho$ throughout, 
we sometimes have to restrict~$\rho$ to subgroups~$H<G$.
In this situation, we say that $V$ is {\it $H$-healthy}\/ 
if the restriction of~$\rho$ to~$H$ is healthy. 
%
\begin{lemma}\label{prep1}
Assume that $H < G$ is a subgroup of finite index. 
Then $V$ is $G$-healthy if and only if $V$ is $H$-healthy.
\end{lemma}

\proof
The implication from $H$-healthy to $G$-healthy is obvious. 
We now assume that $V$ is $G$-healthy and show that $V$ is $H$-healthy as well. 
Since $V$ is $G$-healthy, there exists a vector $v \in V$ and a group element $g \in G$ 
such that $W^g_v$ is infinite-dimensional. 
Denote the right coset~$gH$ in $G/H$ by~$[g]$.

We first consider the special case where the subgroup $H$ is normal. 
Then $G/H$ is a group, and its order $n := |G/H|$ is finite by assumption. 
Hence $[g]^n = id \in G/H$, or equivalently $h := g^n \in H$. 
By Lemma~\ref{linalg} we conclude from the fact that $W^g_v$ is infinite-dimensional 
that $W^h_v$ is infinite-dimensional as well. 
This proves that $V$ is $H$-healthy in the special case that $H \lhd G$ is normal. 

In the general case where $H < G$ is not necessarily normal, 
we consider the normal core of~$H$ in~$G$ defined by
$$
\mathrm{Core}\1 (H) \,:=\, \bigcap_{g \in G} g^{-1} H g .
$$
Note that $\mathrm{Core}\1 (H)$ is a subgroup of~$H$ which is normal in~$G$. 
It is actually the biggest normal subgroup of~$G$ contained in~$H$. 
Moreover, it still has finite index in~$H$, 
see for instance \cite[Theorem~3.3.5]{Sco64}. 
In view of what we already proved, we therefore conclude that
$V$ is $\mathrm{Core}\1 (H)$-healthy. 
Since $\mathrm{Core} \1 (H) < H$ it follows that $V$ is $H$-healthy as well. 
This finishes the proof of the lemma. 
\proofend

If $H < G$ is a subgroup and $v \in V$, we abbreviate by
$$W^H_v=\langle Hv \rangle$$
the subspace of $V$ spanned by the $H$-orbit of the vector $v$. 
The next lemma is our main tool to give an inductive proof of Theorem~\ref{t:healthy}.

\begin{lemma}\label{prep2}
Assume that $H \lhd G$ is a normal subgroup such that $G/H$ is cyclic. 
Suppose further that
$v \in V$ is sick and $W_v^H$ is finite-dimensional. Then $W_v^G$ is finite-dimensional as well. 
\end{lemma}

\proof
Abbreviate $m = \dim W^H_v$. Then there exist $\xi_1, \ldots, \xi_m$ in the group ring $\FF H$ 
such that
\begin{equation} \label{wh}
W^H_v = \big\langle \xi_i v \mid 1 \leq i \leq m \big\rangle.
\end{equation}
Choose $g \in G$ such that $[g] \in G/H$ is a generator.
Assume first that $G/H \cong \ZZ_n$ is finite. Then
$$
W^G_v = \big\langle g^j \xi_i v \mid 1 \leq i \leq m, \, 0 \leq j <n \big\rangle.
$$
Hence $\dim W_v^g \le nm$ is finite.
Assume now that $G/H \cong \ZZ$ is infinite.
The case that $v=0$ is trivial. We therefore assume that~$v \neq 0$. 
Since $v$ is sick by assumption, we have the function $d_v\colon G \to \NN$ from~\eqref{dv}, 
and we abbreviate $n=d_v(g)$. 
Our aim is to show that
\begin{equation}\label{aim}
W^G_v = \big\langle g^j \xi_i v \mid 1 \leq i \leq m, \, 0 \leq j <n \big\rangle.
\end{equation}
For this purpose we abbreviate the right hand side by
$$
X := \big\langle g^j \xi_i v \mid 1 \leq i \leq m, \, 0 \leq j <n \big\rangle.
$$
That $X \subset W_v^G$ is clear. 
We have to check the other inclusion, namely that for every $\eta \in \FF G$ it holds that
$$\eta v \in X.$$
Since the right coset~$[g]$ is a generator of $G/H \cong \ZZ$, 
there exist $\eta_j \in \FF H$ with $\eta_j \neq 0$ 
for only finitely many $j \in \ZZ$ such that
$$
\eta = \sum_{j \in \ZZ} \eta_j g^j.
$$
In view of the definition of $n=d_v(g)$, there exists $\zeta \in \FF G$ of the form
$$
\zeta = \sum_{j=0}^{n-1} \zeta_j g^j, \quad \zeta_j \in \FF H
$$
such that
$$\eta v=\zeta v.$$
Since $H$ is normal in~$G$, there exist $\zeta_j' \in \FF H$ for $0 \leq j <n$ such that
$$
\zeta=\sum_{j=0}^{n-1} g^j \zeta'_j.
$$
In view of \eqref{wh} we conclude that for $0 \leq j <n$ we have
$$
\zeta'_jv \in \big\langle \xi_i v \mid 1 \leq i \leq m \big\rangle
$$
Therefore $\eta v = \zeta v \in X$.
This proves \eqref{aim}. We have shown that
$$
\dim W_v^G \leq nm = d_v(g) \cdot \dim W_v^H,
$$ 
and therefore $W_v^G$ is finite-dimensional. 
The proof of the lemma is complete.
\proofend

We are now in position to prove the main result of this section.

\m \ni
{\it Proof of Theorem~\ref{t:healthy}.}
In view of Lemma~\ref{linalg} it suffices to show that
$V$ is healthy. We argue by contradiction and assume that $V$ is sick. 
By assumption $G$ is virtually polycyclic, hence contains a polycyclic subgroup~$H$ 
of finite index. By Lemma~\ref{prep1} it follows that $V$ is $H$-sick as well. 
Since $H$ is polycyclic, we conclude by applying Lemma~\ref{prep2} inductively that 
$W^H_v$ is finite-dimensional for every $v \in V$. 
Because $H$ has finite index in~$G$ it follows that $W^G_v$ is finite-dimensional for every
$v \in V$ as well. By assumption $V$ is a finitely generated $\FF G$-module. 
Therefore we deduce that $V$ is finite-dimensional, contradicting 
the assumption of the theorem. The proof is complete. 
\proofend

\section{Spicy Hopf algebras} \label{s:spicy}

In this section we consider Hopf algebras~$V$ over a group ring~$\FF G$,
endowed with a filtration~$V^r$, $r >0$,
and exhibit a property of $V$ that guarantees that $\dim V^r$ grows at least like $e^{\sqrt r}$. 

\subsection{Hopf algebras over~$\FF G$}
Let $G$ be a group, $\FF$ a field and $\FF G$ the group ring. 
We first explain the notion of a Hopf algebra over~$\FF G$. 
This might be not completely standard, since usually Hopf algebras are defined over rings 
which are commutative, a requirement that our group ring in general does not fulfill. 
However, the feature which distinguishes a group ring from other non-commutative rings is that 
if $V$ and $W$ are two left modules over~$\FF G$, then we can define on their tensor product
$V \otimes W=V \otimes_{\FF} W$ still the structure of a left $\FF G$-module by using
the tensor of the two representations: 
$g(v \otimes w) := (gv) \otimes (gw)$ for $v \in V$, $w \in W$ and $g \in G$.
A \emph{product}\/ is then an~$\FF G$-linear map $\mu \colon V \otimes V \to V$, 
or equivalently an $\FF$-bilinear map $\mu \colon V \times V \to V$ 
satisfying $\mu(gv,gw) = g \2 \mu(v,w)$ for $g \in G$ and $v,w \in V$.  
Dually, a \emph{coproduct}\/ is then an $\FF G$-linear map $\Delta \colon V \to V \otimes V$. 
To avoid terrible headaches we assume in addition that our product is always associative, 
although this requirement is probably not necessary for the results of this section. 
We abbreviate the product by $vw = \mu(v,w)$. 
If $V$ in addition is graded, i.e.\ $V=\bigoplus_{i=0}^\infty V_i$, where each $V_i$ is an $\FF G$-submodule of~$V$, 
then we grade the tensor product by $(V \otimes V)_k = \bigoplus_{i+j=k}V_i \otimes V_j$ and require in addition that 
the product and coproduct preserve the grading. 
The product endows the tensor product $V \otimes V$ again with a product which is defined
on homogeneous elements by the Koszul sign convention
$$
(v \otimes w)(x \otimes y) \,=\, (-1)^{\deg(w)\deg(x)} vx \otimes wy
$$
where $\deg (v)$ denotes the degree of a homogeneous element~$v$.
Given a left module~$V$ over $\FF G$ and a product~$\mu$ and coproduct~$\Delta$ as above, 
we call the triple $(V,\mu,\Delta)$ a \emph{bialgebra}\/ over~$\FF G$ if $\mu$ and~$\Delta$
are compatible in the sense that 
$\Delta \colon V \to V \otimes V$ is a homomorphism of algebras. 
The bialgebra $(V,\mu,\Delta)$ is called \emph{connected}\/ if $V_0 = \FF$ is one-dimensional 
and if $1 \in \FF$ is also the unit for the multiplication~$\mu$. 
\begin{definition}
{\rm
A connected graded bialgebra $(V,\mu,\Delta)$ over $\FF G$ is called a \emph{Hopf algebra}\/ 
over~$\FF G$ if for every homogeneous element~$v$ of positive degree~$\deg (v)$ the coproduct satisfies
\begin{equation} \label{hopf}
\Delta v \,=\, 1 \otimes v+v \otimes 1 +\sum v_i \otimes v_i'
\end{equation}
with $v_i$ and $v_i'$ of positive degree. 
A vector $v \in V$ is {\it primitive}\/ if $\Delta v = 1 \otimes v+v \otimes 1$.
}
\end{definition}

\subsection{Filtrations} \label{s:filtrations}
To make our Hopf algebra spicy we shall suppose that both the vector space~$V$ and 
the group ring~$\FF G$ are filtered. 
More precisely, we assume that $V$ can be exhausted by a sequence of finite-dimensional 
vector spaces, i.e., for every real number~$r>0$ there exists a finite-dimensional subspace $V^r \subset V$ 
such that $V^r \subset V^s$ for $r \leq s$ and $V = \bigcup_{r>0}V^r$.
Define the {\it value}\/ of $v \in V$ by
$$
|v| := \min \{r \mid v \in V^r \}.
$$
Notice that for scalars $a_i \in \FF$ and vectors $v_i \in V$ we have
\begin{equation} \label{ine:norm}
| a_1 v_1 + \dots + a_n v_n | \,\le\, \max \{ |v_i| \} .
\end{equation}
Dear reader, please do not confuse the value~$|v|$ of~$v$ and its degree~$\deg (v)$, in case $v$ is homogeneous. 
They are not related to each other. 
Also note that the grading is indicated by a subscript 
while the filtration degree is indicated by a superscript. 
To get a filtration on the group ring as well, suppose that the group~$G$ is endowed with a 
{\it length function}, namely a function $L \colon G \to \RR_{\ge 0}$ satisfying
$$
L(g) = L(g^{-1}), \quad L(gh) \leq L(g)+L(h), \qquad g,h \in G.
$$
In the following we abbreviate $|g|=L(g)$. 
Via the length function we can define a filtration on the group ring:
For $r \geq 0$ we define $\FF G^r$ to be the subvector space of~$\FF G$ consisting of finite sums 
$\xi=\sum_{g \in G} \xi_g \1 g$ satisfying $\xi_g=0$ whenever~$|g|>r$.
\begin{definition}
{\rm
The Hopf algebra $(V,\mu,\Delta)$ over~$\FF G$ is called {\it spicy}\/ 
if the vector space~$V$ is endowed with a filtration such that for $v,w \in V$ and $g \in G$ 
it holds that
$$
|vw| \leq |v|+|w|, \qquad |gv| \leq |g|+|v| .
$$
}
\end{definition}
%


We next introduce a condition on spicy Hopf algebras which will guarantee nontrivial lower bounds 
for the growth of $\dim V^r$. 
\begin{definition}
{\rm
Assume that $(V,\mu,\Delta)$ is a spicy Hopf algebra over~$\FF G$. 
A \emph{primitive sequence}\/ is a sequence $(v_i)_{i \in \NN}$
satisfying the following conditions.

\s
\begin{itemize}
\item[(i)] 
The vectors $v_i$, $i \in \NN$, are linearly independent, primitive, 
and of equal positive degree.

\s
\item[(ii)] 
There exists a constant $c>0$ such that $|v_i| \leq c \2 i$.
\end{itemize}
}
\end{definition}

\begin{remark}
{\rm
The notion of a primitive sequence does not involve the action of the group~$G$. 
However, we shall see later that the action of~$G$ is useful to construct primitive sequences. 
}
\end{remark}

\begin{theorem} \label{t:spicy}
Assume that $(V,\mu,\Delta)$ is a spicy Hopf algebra over~$\FF G$ 
which admits a primitive sequence $(v_i)_{i \in \NN}$. 
Then the function $r \mapsto \dim V^r$ grows at least like $e^{\sqrt r}$.

\end{theorem}

\proof
Let $I = \{ i_1, \ldots, i_\ell \} \subset \NN$ be a finite subset of distinct numbers
which we totally order by $i_1<i_2< \ldots<i_\ell$. We abbreviate
$$
v_I = v_{i_1} v_{i_2} \cdots v_{i_\ell} \in V .
$$
Let $m$ be the common degree of the vectors~$v_i$.
In view of property~(i) of a primitive sequence, and since $\deg {v_{I}} = m^{\#I}$,
Proposition~\ref{p:quantpbw} below shows that 
the vectors $v_I$, $I \subset \NN$, are linearly independent.
We can assume without loss of generality that the constant~$c$ in property~(ii) 
of a primitive sequence is~$1$. 
Then we have
$$
|v_I| \leq \sum_{j=1}^\ell |v_{i_j}| \leq \sum_{j=1}^\ell i_j.
$$
For $n \in \NN$ denote by $q(n)$ the number of partitions of $n$ into distinct parts. We have shown that
$$
\dim V^n \geq q(n).
$$
By Euler's theorem, the number of partitions of~$n$ into distinct parts coincides with 
the number of its partitions into odd parts,
see for example \cite[Corollary 1.2]{andrews}. 
The asymptotics of this sequence coincides up to a constant with the asymptotics of 
the partition function (see for example \cite[Chapter 16]{nathanson}),
which grows like $e^{C \sqrt{n}}$ for a positive constant~$C$
according to a theorem of Hardy and Ramanujan, see for example \cite[Chapter 15]{nathanson}.  
\proofend


Combining Theorems~\ref{t:healthy} and~\ref{t:spicy} we obtain the following result.

\begin{corollary} \label{c:growth}
Assume that $(V,\mu,\Delta)$ is a spicy Hopf algebra over~$\FF G$ 
where $G$ is a virtually polycyclic group.
Assume further that $\oplus_{i < m} V_i$ is finite-dimensional and that 
$V_m$ is infinite-dimensional but finitely generated as an $\FF G$-module.
Then the function $r \mapsto \dim V^r$ grows at least like $e^{\sqrt r}$.
\end{corollary}

\proof
By Theorem~\ref{t:healthy} there exists $v \in V_m$ and $g \in G$ such that
the vectors $v_i := g^i v$, $i \in \NN$, are linearly independent.
The vectors $v_i$ have equal degree~$m$.
Using the defining properties of a filtration, we estimate
$$
|v_i| \,=\, |g^iv| \,\le\, |g^i| + |v| \,\le\, i \2 |g| + |v|.
$$
With $c := |g| + |v|$ we thus have $|v_i| \le c \2 i$.
The sequence $(v_i)_{i \in \NN}$ therefore meets all the properties of a primitive sequence, 
except that the $v_i$ may fail to be primitive. 
To correct this, we consider the linear map 
$A \colon \langle v_1, v_2, \dots \rangle \to  V_m \otimes V_m$ given by
$$
A(v) = \Delta v - 1 \otimes v - v \otimes 1 .
$$
Since the vector space $\oplus_{j < m} V_j$ is finite-dimensional, 
the subvector space of $V_m \otimes V_m$ spanned by the elements 
$u_i \otimes u_i'$ with $u_i, u_i' \in \oplus_{j < m}V_j$ is finite-dimensional.
Since $A$ takes values in this finite-dimensional vector space,
$k := \rank A$, which equals the dimension of $\langle v_1, v_2, \dots \rangle / \ker A$, is finite.
In particular, $\ker A$ in infinite-dimensional.
We shall construct by induction a sequence $w_1, w_2, \dots$ of linearly independent 
elements in $\ker A$ such that $|w_i| \le c (k+1)$. 
The sequence $(w_i)_{i \in \NN}$ is then a primitive sequence in~$V$, and the corollary
follows in view of Theorem~\ref{t:spicy}.

The restriction of the map $A$ to $\langle v_1, v_2, \dots, v_{k+1} \rangle$ has a non-trivial kernel.
Let $w_1 := a_1 v_1 + \dots + a_{k+1} v_{k+1}$ be a non-trivial element in this kernel.
Then by \eqref{ine:norm}, $|w_1| \le c(k+1)$.
Next, the restriction of $A$ to $\langle v_{k+2}, v_{k+3}, \dots, v_{2(k+1)} \rangle$ has a non-trivial kernel.
Let $w_2 := a_{k+2} v_{k+2} + \dots + a_{2(k+1)} v_{2(k+1)}$ be a non-trivial element in this kernel.
Then by \eqref{ine:norm}, $|w_2| \le 2c(k+1)$, and $w_1, w_2$ are linearly independent
because $v_1, \dots, v_{2(k+1)}$ are linearly independent.
Proceeding in this way we construct linearly independent vectors $w_1, w_2,\dots$ such that
$|w_i| \le i \1 c \1 (k+1)$.
\proofend


\section{A quantum Poincar\'e--Birkhoff--Witt theorem} \label{s:PBW}

\ni
Consider a graded bialgebra $(V,\Delta,\mu)$ over the field~$\FF$ 
which is connected (i.e., $V_0=\FF$ is one-dimensional) and is such that 
$1 \in \FF$ also serves as the unit for the multiplication~$\mu$.
(The group~$G$ plays no role in this section.) 
We again assume that $\mu$ is associative, and write $vw=\mu(v,w)$. 
Also recall that $v \in V$ is \emph{primitive}\/ if
$\Delta(v)=1 \otimes v+v \otimes 1$. 
Suppose that for $N \in \NN$ we are given linearly independent and primitive 
elements $v_1, \ldots ,v_N \in V_m$ of equal positive degree $m \in \NN$. 
Abbreviate $\NN_N=\{1, \ldots, N\}$.
We order~$I  \subset \NN_N$ using the canonical order of~$\NN_N$, 
namely we write $I=\{ i_1, \ldots, i_\ell \}$ satisfying $i_1<i_2< \ldots <i_\ell$. 
We then abbreviate $v_I = v_{i_1} v_{i_2} \cdots v_{i_\ell}$, where we use the convention that $v_\emptyset=1$. 
The following proposition reminiscent of the Poincar\'e--Birkhoff--Witt theorem 
seems to be known to people working in quantum group theory, 
cf.\ \cite[Theorem~1.5(b)]{wang-zhang-zhuang}. 
However, since the Hopf algebras arising in the theory of quantum groups are usually
not graded, we provide a proof for the readers convenience.

\begin{proposition} \label{p:quantpbw}
Assume that $v_1, \dots, v_N$ are linearly independent and primitive vectors in~$V$.
Then the vectors $v_I$, $I \subset \NN_N$, are linearly independent.
\end{proposition}

\proof
The crucial ingredient in the proof is the computation of the coproduct of
the elements~$v_I$. To determine the signs in this formula, we use the following convention. 
Given a subset $I$ of $\{1, \ldots, N\}$ we order $I=\{i_1, \ldots, i_\ell\}$ and its complement 
$I^c=\{j_1, \ldots, j_{N-\ell}\}$. 
This determines a permutation $(1, \ldots, N) \mapsto (i_1, \ldots, i_\ell, j_1, \ldots, j_{N-\ell})$. 
We denote by $\sigma(I)$ the signum of this permutation. 
\begin{lemma} \label{le:I}
The coproduct of $v_{\NN_N}$ is given by
$$
\Delta(v_{\NN_N}) =\sum_{I \subset \NN_N} \sigma(I)^m \2 v_I \otimes v_{I^c}.
$$ 
Hence if $m$ is even, then
$\Delta(v_{\NN_N}) =\sum_{I \subset \NN_N} v_I \otimes v_{I^c}$, 
and if $m$ is odd, then
$\Delta(v_{\NN_N}) = \sum_{I \subset \NN_N} \sigma(I) v_I \otimes v_{I^c}$.
\end{lemma}

\proof 
We prove the lemma by induction on~$N$. 
For $N=1$ the lemma is an immediate consequence of the fact that 
$v_1$ is primitive. 
For the induction step we assume that the formula holds for~$N-1$. 
We compute 
\begin{eqnarray*}
\Delta(v_{\NN_N})&=&\Delta(v_{\NN_{N-1}}v_N) \\
&=&\Delta(v_{\NN_{N-1}}) \2 \Delta(v_N) \\
&=&\Bigg(\sum_{I \subset \NN_{N-1}} \sigma(I)^m \2 v_I \otimes v_{\NN_{N-1} \setminus I}\Bigg)\Bigg(1 \otimes v_N+v_N \otimes 1\Bigg)\\
&=&\sum_{I \subset \NN_{N-1}} \sigma(I)^m \2 v_I \otimes (v_{\NN_{N-1} \setminus I} \2 v_N)\\
& &+\sum_{I \subset \NN_{N-1}} (-1)^{\deg (v_{\NN_{N-1}\setminus I}) \deg (v_N)}\sigma(I)^m \2 (v_I v_N) \otimes v_{\NN_{N-1} \setminus I} .
\end{eqnarray*}
Since $\deg (v_N) =m$ and $\deg (v_{\NN_{N-1}\setminus I}) = m (N-1-\# I)$, 
we have for each $I \subset \NN_{N-1}$ that
\begin{eqnarray*}
(-1)^{\deg (v_{\NN_{N-1}\setminus I}) \deg (v_N)} \2 \sigma(I)^m &=& (-1)^{m^2(N-1-\#I)} \2 \sigma(I)^m \\ 
&=& (-1)^{m(N-1-\#I)} \2 \sigma(I)^m \\ 
&=& \sigma(I \cup \{N\})^m .
\end{eqnarray*}
The previous sum therefore becomes
\begin{eqnarray*}
&=&\sum_{I \subset \NN_{N-1}} \sigma(I)^m \2 v_I \otimes v_{\NN_N \setminus I}\\
& &+\sum_{I \subset \NN_{N-1}} \sigma(I \cup \{N\})^m \2 v_{I \cup \{N\}} \otimes v_{\NN_N \setminus (I \cup \{N\})}\\
&=&\sum_{I \subset \NN_N} \sigma(I)^m \2 v_I \otimes v_{\NN_N \setminus I}.
\end{eqnarray*}
This proves the induction step and hence the lemma. 
\proofend

\ni
{\it Proof of Proposition~\ref{p:quantpbw}.} 
Recall that by assumption the vectors $v_1, \dots, v_N$ are linearly independent and have all the same
degree~$m$.
By looking at the degree we see that it suffices to show that for every $k \leq N$ 
the vectors~$v_I$ with $I \subset \NN_N$ and $\# I = k$ are linearly independent. 
For $k=1$ this is an assumption. 
For the induction step we assume this assertion for all $k \leq n-1$ where $n \leq N$. 
Since the coproduct is linear, it suffices to show that the vectors $\Delta (v_I)$ 
with $I \subset \NN_N$, $\#I=n$ are linearly independent. 
It follows from the induction hypothesis that for every $0<k<n$ the vectors 
$v_I \otimes v_J$ with $I,J \subset \NN_N$ and $\#I=k$, $\#J=n-k$ are linearly independent. 
This and the formulas
$$
\Delta(v_I) = \sum_{J \subset I} \sigma(J)^m \2 v_J \otimes v_{I \setminus J},
$$
that follow from (the proof of) Lemma~\ref{le:I},
complete the induction step. 
\proofend

\section{Proof of Theorem~\ref{t:main}} \label{s:proof}

\ni
We shall deduce Theorem~\ref{t:main} from Corollary~\ref{c:growth}.
We start with saying briefly ``who is who'' in Corollary~\ref{c:growth}:
Let $M$ be a closed connected manifold, and let $\widetilde M$ be its universal covering space.
We take as $G$ the fundamental group $\pi_1(M)$ of~$M$,
and as $V$ we take the homology $H_*(\Omega \widetilde M; \FF)$ over a suitable field $\FF$.
The product~$\mu$ will be the Pontryagin product given by concatenation of loops, 
and the coproduct will simply come from the diagonal map 
$\Omega \widetilde M \to \Omega \widetilde M \times \Omega \widetilde M$, 
$x \mapsto (x,x)$.
Now fix a Riemannian metric on~$M$. 
The filtration on $V^r$ will be given by taking $|v|$ as the smallest~$r$ for which $v$
can be represented by a cycle of based loops of length at most~$r$,
and for $g \in \pi_1 (M)$ we take $L(g)$ to be (half of) the length of 
the shortest curve representing~$g$.

We shall show in Sections~\ref{s:Hopf} and~\ref{s:filtration} that with these choices,
$(V,\mu,\Delta)$ is a spicy Hopf algebra over~$\FF G$.
In Section~\ref{s:dim} we show that for $M$ of non-finite type with $G=\pi_1(M)$ 
virtually polycyclic, 
the dimension assumptions on $\oplus V_i$ are met.

\subsection{The Hopf-algebra structure on $H_* (\Omega \widetilde M;\FF)$}  \label{s:Hopf}

Let $M$ be a closed connected manifold, and let $\pr \colon \widetilde M \to M$ be its universal covering space.
Fix $p \in M$ and $\tilde p \in \widetilde M$ over~$p$.
The spaces $\Omega_0M$ and $\Omega \widetilde M$ of contractible continuous loops based at~$p$
and~$\tilde p$, respectively, are canonically identified.
Conjugation of loops in $\Omega_0M$ by loops in~$M$ based at~$p$ yields an action of 
the fundamental group $G = \pi_1(M,p)$ on $H_*(\Omega_0M) = H_* (\Omega \widetilde M)$:
Given a cycle $C = \{\gamma\}$ of loops in $\Omega \widetilde M$ based at~$\tilde p$
and given $g \in G$,
the class $g [C]$ is defined as the class represented by the cycle of loops $\{c_g^{-1} \circ g\gamma \circ c_g\}$,
where $c_g$ is the lift to~$\widetilde M$ starting at~$\tilde p$ of a loop in~$M$ in class~$g$,
and $g \gamma$ is the lift of $\pr \circ \gamma$ starting at $g \tilde p$.

Let $\FF$ be the field, and abbreviate $V_* := H_* (\Omega \widetilde M;\FF)$.
The action of~$G$ on~$V$ extends to an action of $\FF G$ on~$V$.
Concatenation of loops in $\widetilde M$ based at~$\tilde p$ induces a product 
$\mu \colon V \otimes V \to V$, called the Pontryagin product. 
Since concatenation of loops is associative up to homotopy,
$\mu$ is associative. 
It follows from the definition of the action of~$G$ that 
$\mu$ is $\FF G$-linear.
In order to define the coproduct~$\Delta \colon V \to V \otimes V$,
we consider, more generally, a topological space~$X$ and
the diagonal map $\delta_X \colon X \to X \times X$, $x \mapsto (x,x)$.
Since we work over a field~$\FF$, the cross product 
$$
H_*(X;\FF) \otimes H_*(X;\FF) \,\stackrel{\times}{\longrightarrow}\, H_*(X \times X;\FF)
$$
is an isomorphism by the K\"unneth formula.
We can therefore define $\Delta_X \colon H_*(X;\FF) \to H_*(X;\FF) \otimes H_*(X;\FF)$ by
$$
\Delta_X := \times^{-1} \circ (\delta_X)_* .
$$
Assume now in addition that $X$ is a path-connected H-space, with product~$\nu$.

\begin{lemma} \label{le:Hopf}
The homology $H_*(X;\FF)$, with product induced by~$\nu$ and with coproduct~$\Delta_X$,
is a Hopf algebra.
\end{lemma}

We refer to \cite[Theorem~7.15]{Whi78} for the proof.
For the readers convenience, we verify that for every homogeneous element~$v$ of positive degree, 
$\Delta_X v$ has the form~\eqref{hopf}.
Let $p \colon X \times X \to X$, $(x,y) \mapsto x$, be the projection on the first factor. 
Then
\begin{equation} \label{e:pd}
p \circ \delta_X = id_X .
\end{equation}
For an element $u = v_n' \otimes 1 + 1 \otimes v_n'' + \sum v_i' \otimes v_j'' 
\in \oplus_{i+j=n} H_i(X) \otimes H_j(X)$ 
with $\deg v_i' < n$ we have
\begin{equation} \label{e:px}
(p_n \circ \times ) \1 u = v_n' 
\end{equation}
by the geometric definition of the cross product (see e.g.\ \cite[\S 3.B]{Hat.AT}).
Now write $\Delta_X v = v_n' \otimes 1 + 1 \otimes v_n'' + \sum v_i' \otimes v_j''$
with $\deg v_i', \deg v_j'' \ge 1$.
Using \eqref{e:pd}, the definition of $\Delta_X$ and \eqref{e:px} we get
\begin{eqnarray*}
v &=& p_n \circ (\delta_X)_n v \\
  &=& p_n \circ \times \circ \Delta_X v \\
  &=& p_n \circ \times \left( v_n' \otimes 1 + 1 \otimes v_n'' + \sum v_i' \otimes v_j'' \right) \\
  &=& v_n'.
\end{eqnarray*}
Similarly we find $v_n'' = v$, and so $\Delta_X v = v \otimes 1 + 1 \otimes v + \sum v_i' \otimes v_j''$
with $\deg v_i', \deg v_i'' \ge 1$.
\proofend

Since $\widetilde M$ is simply connected, $\Omega \widetilde M$ is path-connected.
Hence $V_0 = H_0 (\Omega \widetilde M; \FF) \cong \FF$ is one-dimensional.
Moreover, $1 \in \FF$ corresponds to the class of the constant path $\tilde p \in \widetilde M$,
which is the unit for the Pontryagin product~$\mu$.
Applying Lemma~\ref{le:Hopf} with $X = \Omega \widetilde M$ and writing $\Delta$ for $\Delta_{\Omega \widetilde M}$ 
we obtain that $(V,\mu, \Delta)$ is a Hopf algebra over~$\FF G$.
(We have already noticed that $\mu$ is $\FF G$-linear.
The definition of the $G$-actions on~$V$ and on $V \otimes V$ shows that also $\Delta$ is $\FF G$-linear.)

\subsection{The filtration on $H_* (\Omega \widetilde M;\FF)$} \label{s:filtration}
Fix a Riemannian metric~$\rho$ on~$M$, and let $\tilde \rho = \pr^* \rho$ be the corresponding 
Riemannian metric on~$\widetilde M$.
Given a piecewise smooth curve $\gamma$ in $\widetilde M$, we denote by $\ell (\gamma)$
the length of $\gamma$ with respect to the Riemannian metric $\widetilde \rho$.
For $r>0$ let $V^r$ be the set of homology classes in $V = H_*(\Omega \widetilde M;\FF)$
that can be represented by cycles formed by piecewise smooth loops $\gamma$ based at~$\tilde p$ 
with $\ell (\gg) \le r$,
$$
V^r \,:=\, \iota^r_* H_*(\Omega^r \widetilde M; \FF) .
$$
Then each $V^r$ is finite-dimensional (see \cite[\S 16]{Mil63}), 
$V^r \subset V^s$ for $r \le s$ and $V = \bigcup_{r>0}V^r$.
As in Section~\ref{s:filtrations} define the value of $v \in V$ by 
$|v| := \min \left\{ r \mid v \in V^r \right\}$.
In view of the definition of the Pontryagin product and by the triangle inequality, 
$|v w| \le |v| + |w|$ for all $v,w \in V$.

Next, for $g \in G = \pi_1(M,p)$ let $\ell (g)$ be the minimal length of a piecewise smooth loop based at~$p$
that represents~$g$. 
In other words, $\ell (g)$ is the length of the shortest geodesic lasso based at $p$ in class~$g$.
Set $L(g) := \frac 12 \ell (g)$. 
Then $L(g) = L(g^{-1})$ and $L(gh) \le L(g) + L(h)$ for all $g,h \in G$ by the triangle inequality.
Finally, for $g \in G$ denote by $c_g$ the lift to $\widetilde M$ based at $\tilde p$ of a shortest curve 
in class~$g$.
Then
$$
\ell (c_g^{-1} \circ g \gamma \circ c_g) \,\le\,
\ell (c_g^{-1}) + \ell (g \gamma) + \ell (c_g) \,=\, \ell (\gamma) + 2 \ell (c_g) \,=\, \ell (\gamma) + L(g)
$$
for all $g \in G$ and $\gamma \in \Omega \widetilde M$.
In view of the definition of the $G$-action on~$\Omega \widetilde M$ we find that
$|g v| \le |v| + L(g) = |v| + |g|$.
We have shown that $(V,\mu, \Delta)$ is a spicy Hopf algebra over~$\FF G$.

\subsection{Dimensions} \label{s:dim}

Recall that $M$ is of {\it finite type}\/ if its universal cover~$\widetilde M$
is homotopy equivalent to a finite CW-complex.

\begin{lemma} \label{le:finitetype}
{\rm (\cite[Lemma~2.2]{FraLabSch13})}
The following are equivalent.

\begin{enumerate}
\item[(i)]
$M$ is of finite type.

\s
\item[(ii)]
The Abelian groups $H_k(\widetilde M)$ are finitely generated for all $k \ge 1$.

\s
\item[(iii)]
The Abelian groups $\pi_k(M)$ are finitely generated for all $k \ge 2$.
\end{enumerate}
\end{lemma}

Now assume that $M$ is not of finite type.
By the lemma, we can define 
\begin{equation} \label{def:m}
\fm (M) \,:=\, 
\min \left\{ k \mid H_k (\widetilde M) \mbox{ is not finitely generated} \right\} \,\in\, \{ 2, \dots, \dim M \} .
\end{equation}

The main result of this subsection is

\begin{proposition} \label{p:HOmega}
Assume that $M$ is not of finite type. 
Let $\fm = \fm(M)$ be as in definition~\eqref{def:m}. Then

\begin{itemize}
\item[(i)]
$H_i (\Omega \widetilde M)$ is finitely generated for $i \le \fm-2$;

\s
\item[(ii)]
$H_{\fm-1} (\Omega \widetilde M)$ is not finitely generated,
but is finitely generated as a $\ZZ G$-module.
\end{itemize}
\end{proposition}

\proof
For each $k \ge 2$ the fundamental group $G = \pi_1(M,p)$ acts on~$\pi_k(M) = \pi_k(M,p)$ by conjugation.
Under the Hurewicz homomorphism $h \colon \pi_k(M) = \pi_k(\widetilde M) \to H_k(\widetilde M)$
this action corresponds to the action of $G$ on~$H_k(\widetilde M)$ induced by
deck transformations (the commutative diagram on the left).
On $\pi_{k-1}(\Omega \widetilde M) = \pi_{k-1}(\Omega_0M) \cong \pi_k(M)$ 
this action of~$G$ is induced by conjugation of elements in $\Omega_0M \cong \Omega \widetilde M$. 
This action also induces an action on $H_{k-1}(\Omega \widetilde M)$, 
namely the action described in Section~\ref{s:Hopf},
and the two actions commute with the Hurewicz homomorphism 
$h \colon \pi_{k-1}(\Omega \widetilde M) \to H_{k-1}(\Omega \widetilde M)$
(the commutative diagram on the right),
$$
\xymatrixcolsep{3pc}\xymatrix{ 
H_k(\widetilde M)
\ar[d]  &
\pi_k(\widetilde M) = \pi_{k-1}(\Omega \widetilde M) \ar[d] \ar[l]_-{\;\;\;h} \ar[r]^-h &
H_{k-1}(\Omega \widetilde M) \ar[d] 
\\
H_k(\widetilde M)  & 
\pi_k(\widetilde M) = \pi_{k-1}(\Omega \widetilde M) \ar[l]^-{\;\;\;h} \ar[r]_-h &
H_{k-1}(\Omega \widetilde M)
}
$$

Recall that $\fm \in \{ 2, \dots, \dim M \}$ is the minimal integer such that 
$H_\fm(\widetilde M)$ is not finitely generated.
By Serre's theory of $\cc$-classes, applied to the class of finitely generated Abelian groups, 
$\fm$ is also the minimal integer such that $\pi_\fm(\widetilde M)$ is not finitely generated.
More precisely, Serre's Hurewicz theorem implies that for $k \le \fm$ the Hurewicz map 
$h \colon \pi_k(\widetilde M) \to H_k(\widetilde M)$ is injective and surjective up to finitely generated groups, 
see \cite{Ser53} or \cite[Theorem~15 on p.~508]{Spa66} or \cite[Theorem~1.8]{Hat.spec}.
Hence $\pi_k (\Omega \widetilde M)$ is finitely generated for $k \le \fm-2$, 
but not so for $k = \fm-1$.

Since $\Omega \widetilde M$ is a path-connected H-space, 
$\pi_1 (\Omega \widetilde M)$ acts trivially on $\pi_k (\Omega \widetilde M)$ 
for $k \ge 0$.  
Serre's Hurewicz theorem now implies that $h \colon \pi_k(\Omega \widetilde M) \to H_k(\Omega \widetilde M)$ 
has finitely generated kernel and cokernel for $k \le \fm-1$,
see \cite[p.~274]{Ser53} or \cite[Theorem~20 on p.~510]{Spa66} or \cite[Theorem~1.8]{Hat.spec}.
It follows that $H_k(\Omega \widetilde M)$ is finitely generated for $k \le \fm-2$, 
but not so for $k = \fm-1$.

\s
We are left with proving the second assertion in~(ii).
After replacing~$M$ by a homotopy equivalent space, if necessary, 
we find a CW-structure on~$M$. Since $M$ is compact, this CW-structure is finite.
We lift this structure to~$\widetilde M$ by the action of~$G$.
The cellular chain complex $C_*(\widetilde M;\ZZ)$ is then a finitely generated $\ZZ G$-module 
in each degree.
Since $G$ is virtually polycyclic, 
the ring~$\ZZ G$ is left Noetherian,
see~\cite{Hal54} or \cite{Iva90}. 
Hence each $\ZZ G$-module $C_*(\widetilde M;\ZZ)$ is left Noetherian.
Therefore the kernel and the image of the differential of $C_*(\widetilde M;\ZZ)$ 
as well as the quotient $H_*(\widetilde M)$
are finitely generated left Noetherian $\ZZ G$-modules in each degree.
In particular, $H_{\fm}(\widetilde M)$ is a finitely generated $\ZZ G$-module.
Recall that $h \colon \pi_{\fm}(\widetilde M) \to H_{\fm}(\widetilde M)$ 
is injective and surjective up to finitely generated groups.
In view of the commutative diagram above, 
it follows that $\pi_{\fm}(\widetilde M)$ and hence $\pi_{\fm-1}(\Omega \widetilde M)$ are finitely generated $\ZZ G$-modules.
As we have seen before, $h \colon \pi_{\fm-1}(\Omega \widetilde M) \to H_{\fm-1}(\Omega \widetilde M)$ 
has finitely generated cokernel. 
Hence $H_{\fm-1} (\Omega \widetilde M)$ is also a finitely generated $\ZZ G$-module.
\proofend

\subsection{Proof of Theorem~\ref{t:main} and Remark~\ref{rem:HP}.1}  \label{s:end}
By the universal coefficient theorem, $H_i (\Omega \widetilde M; \FF) = H_i (\Omega \widetilde M) \otimes \FF$
for every field~$\FF$.
Hence assertion~(i) implies that $H_i (\Omega \widetilde M;\FF)$ is finite-dimensional for $i \le \fm-2$,
and assertion~(ii) implies that $H_{\fm-1} (\Omega \widetilde M;\FF)$ is finitely generated as an $\FF G$-module. 
Assume now that $\FF$ is a field such that $H_{\fm} (\widetilde M;\FF) = H_{\fm}(\widetilde M) \otimes \FF$ is infinite-dimensional.
As we have seen in the proof above, the two Hurewicz maps 
\begin{equation} \label{e:hh}
\xymatrixcolsep{3pc}\xymatrix{ 
H_{\fm}(\widetilde M) &
\pi_{\fm}(\widetilde M) = \pi_{\fm-1}(\Omega \widetilde M) \ar[l]_-{\;\;\;h} \ar[r]^-h &
H_{\fm-1}(\Omega \widetilde M) 
}
\end{equation}
both have finitely generated kernel and cokernel.
If follows that $H_{\fm-1}(\Omega \widetilde M) \otimes \FF = H_{\fm-1} (\Omega \widetilde M;\FF)$ 
is also infinite-dimensional.
Hence the dimension assumptions in Corollary~\ref{c:growth} are satisfied, and we conclude
that $\dim V^T$ grows at least like $e^{\sqrt T}$.

Similarly, since the left map in~\eqref{e:hh} has finitely generated kernel and cokernel,
$H_{\fm}(\widetilde M;\FF) = H_{\fm}(\widetilde M) \otimes \FF$ is infinite-dimensional if and only if 
$\pi_{\fm}(\widetilde M) \otimes \FF$ is infinite-dimensional.
\proofend

\section{Examples} \label{s:exa}

\subsection{A class of examples} \label{ss:exa.neu}
In this paragraph we construct examples of manifolds with infinite cyclic fundamental group
that meet the assumptions of Theorem~\ref{t:main}.

Consider a simply-connected manifold~$X$ of dimension $d \ge 4$
that is not homeomorphic to a sphere.
Let $m \ge 2$ be the minimal $k$ such that $H_k(X) \neq 0$.
Then $m \le d/2 \le d-2$ by Poincar\'e duality.
Let $\check X$ be the compact manifold with boundary obtained by removing from~$X$
the interior of two disjoint closed $d$-balls $B_1 \cup B_2 \subset X$.
Then
\begin{equation}  \label{e:HX}
H_k(\check X) \,=\, H_k(X) \quad \mbox{ for } \, 1 \le k \le m
\end{equation}
since $H_k(\pp B_i) = H_k(S^{d-1}) =0$ for these~$k$.
Choose a diffeomorphism $\gf \colon \pp B_1 \to \pp B_2$,
and let $M$ be the manifold obtained from~$\check X$ by identifying $\pp B_1$ with~$\pp B_2$ via~$\gf$.
Then $\pi_1(M) = \ZZ$ by the Seifert--van Kampen theorem.
For $n \in \ZZ$ let $\check X_n$ be a copy of~$\check X$.
The universal cover~$\widetilde M$ is obtained by glueing ``the right boundary'' $\pp B_1$ of~$\check X_n$ by~$\gf$
to ``the left boundary'' $\pp B_2$ of~$\check X_{n+1}$ for $n \in \ZZ$.
For $N \in \NN \cup \{0\}$ consider the part $\widetilde M_N := \bigcup_{-N \le n \le N} \check X_n$ of~$\widetilde M$.
By the Mayer--Vietoris theorem and by~\eqref{e:HX},
$$
H_k(\widetilde M_N) \,=\, \oplus_{n=-N}^N H_k (\check X_n) \,=\, \oplus_{2N+1} H_k(X) \quad \mbox{ for } \, 1 \le k \le m. 
$$
Hence $H_k(\widetilde M) = \varinjlim H_k(\widetilde M_N) = \oplus_{\ZZ} H_k(X)$ for $1 \le k \le m$.
Recalling that $H_k(M) = 0$ for $1 \le k \le m-1$ and $H_m(M) \neq 0$, we see that 
$M$ is not of finite type, and that $m = \fm = \fm(M)$.
For every field~$\FF$ we have
$$
H_\fm(\widetilde M;\FF) \,=\, \oplus_{\ZZ} H_\fm(X) \otimes \FF .
$$
Moreover, $H_\fm(X)$ is a non-trivial finitely generated Abelian group,
$$
H_\fm(X) \,\cong\, \ZZ^r \oplus \ZZ_{q_1} \oplus \dots \oplus \ZZ_{q_\ell}
$$
with $r \ge 0$ and the $q_i$ powers of primes.
If $r \ge 1$ choose $\FF = \QQ$. Then $H_\fm (\widetilde M;\QQ) = \oplus_\ZZ  H_\fm(X) \otimes \QQ = \oplus_\ZZ \QQ^r$
is infinite-dimensional.
If $r=0$ let $p$ be the prime number dividing~$q_1$.
Then $H_\fm (\widetilde M;\FF_p)$ contains $\oplus_\ZZ \FF_p$ as a subvector space, and hence is also infinite-dimensional.

We see that the assumptions of Theorem~\ref{t:main} are met, and conclude that for a suitable field~$\FF$ the dimension of
$\iota_* H_* (\Omega_0^T M;\FF)$
grows at least like $e^{\sqrt T}$.
While this seems to be a new result if $X$ or, equivalently, $M$ is irreducible, 
it has been shown in the proof of Theorem~D in~\cite{PatPet04} 
that if~$M$ can be written in the form $M = M_1 \# M_2$ with $M_2$ simply connected and not homeomorphic to a sphere, then 
$\dim \iota_* H_* (\Omega_0^T M;\QQ)$ grows even exponentially.

\subsection{A ``counterexample''} \label{ss:exa}
We illustrate the role of our standing assumption that there exists a field~$\FF$ such that 
$H_{\fm}(\widetilde M;\FF)$ is infinite-dimensional by an example
(in which $M$ is not a manifold, but a CW-complex).
Let $M$ be the mapping torus of a degree-two map~$f$ of the 2-sphere $S^2$.
Then $\pi_1(M) = \ZZ$, and $\widetilde M$ is the ``double mapping telescope'' obtained 
by glueing together the mapping cylinders of $f_i = f$, $i \in \ZZ$.
This space deformation retracts onto the mapping telescope formed by the mapping cylinders 
of $f_i = f$, $i \ge 0$. Therefore $H_2(\widetilde M)$ can be identified with $\ZZ[1/2]$, 
the subgroup of~$\QQ$ consisting of rational numbers with denominators a power of~$2$
(see Exercise~1 of Section~3.F and Example~3.F.3 in~\cite{Hat.AT}).
The Abelian group $\ZZ[1/2]$ is not finitely generated, hence $\fm =2$ in this example. 
However, for a field~$\FF$ of characteristic~$p$ we have $H_2(\widetilde M;\FF) =0$ if $p=2$ 
and $H_2(\widetilde M;\FF) \cong \FF$ otherwise.

The group $\ZZ [1/2]$ has the property that
every finitely generated subgroup is generated by one element 
(indeed, the smallest positive element in the subgroup is a generator).
In particular, $\ZZ[1/2]$ is the union of a nested sequence of subgroups generated by
one element (for instance the subgroups generated by $1/2^j$, $j \ge 1$).
These properties strongly distinguish $\ZZ[1/2]$ from an infinite-dimensional vector space.

We have seen in the proof of Proposition~\ref{p:HOmega} that $V_1 := H_1(\Omega \widetilde M)$ 
is finitely generated as a $\ZZ[\pi_1 M]$-module. 
This can be seen very explicitely in this example:
Since $\pi_1(\Omega \widetilde M) \cong \pi_2(\widetilde M)$ is Abelian, 
$V_1 = H_1(\Omega \widetilde M) = \pi_1(\Omega \widetilde M) = \pi_2(\widetilde M) = H_2(\widetilde M) \cong \ZZ[1/2]$.
Let $t = 1$ be the generator of the fundamental group $G = \ZZ$ of~$M$.
The action of $t \in G$ on $H_2(\widetilde M)$ and hence on~$V_1$ corresponds to multiplication by~$2$
on~$\ZZ[1/2]$. 
Moreover, 
the group ring~$\ZZ G$ is the ring of Laurent polynomials $\ZZ[t,t^{-1}]$.
It follows that $V_1 \cong \ZZ[1/2]$ is generated by~$1$ over~$\ZZ G$.
We see that ``the $\ZZ$-version'' of the assumptions in Theorem~\ref{t:healthy} 
is satisfied.
While it is still true that for every non-zero vector $v \in V_1$ the subgroup 
generated by $\{ t^i v \}_{i \in \ZZ}$ is not finitely generated, 
we cannot use the sequence $(t^i v)_{i \in \ZZ}$ to prove a lower bound on the rank of~$V^r$
by the arguments of Sections~\ref{s:spicy} and~\ref{s:PBW},
because for each~$k$ the subgroup generated by $\{t^i v\}_{|i| \le k}$
has only rank one (being generated by $t^{-k}v$).


\section{Proof of Corollaries~\ref{c:chords} and~\ref{c:ent}} \label{s:cor}

\m \ni
{\it Proof of Corollary~\ref{c:chords}.}
Let $\cc_{pq}$ be a component of $\Omega_{pq}M$.
Fix a smooth path $c \in \cc_{pq}$, of length $\ell (c)$. 
The map $h \colon \Omega_0(M,p) \to \cc_{pq}$, $\gg \mapsto c \circ \gg$ is a homotopy equivalence.
It maps $\Omega_0^TM$ to $\cc_{pq}^{T+\ell(c)}$, where $\cc_{pq}^T$ is the space of piecewise smooth paths in~$\cc_{pq}$
of length~$\le T$.
Theorem~\ref{t:main} now implies that 
\begin{equation} \label{e:cce}
\dim \iota_*^T H_* (\cc_{pq}^T;\FF) \quad \mbox{grows at least like $e^{\sqrt T}$.}
\end{equation}
Notice that $\dim \iota_*^T H_* (\cc_{pq}^T;\FF) \le \dim H_*(\cc_{pq}^T;\FF)$.
For geodesic flows, Corollary~\ref{c:chords} now follows from classical Morse theory, 
see \cite[Theorem~16.3]{Mil63} or \cite[p.\ 116]{Pat:book}.
For the general case of Reeb flows, we use~\eqref{e:cce} and Lagrangian Floer homology,
exactly as in \cite[Section~6]{MacSch11}.
\proofend


\m \ni
{\it Proof of Corollary~\ref{c:ent}.}
For geodesic flows the claim follows from Corollary~\ref{c:chords} 
and the geometric arguments in~\cite[Section~3.1]{Pat:book}.
For the general case of Reeb flows we use Theorem~\ref{t:main} and
the idea from~\cite{FS:GAFA}, 
that was further developed in~\cite{FraSch06,MacSch11,FraLabSch13}.
The proof goes exactly as the proof of Theorem~4.6 in~\cite{MacSch11};
we therefore only sketch the proof.
Let $\Sigma \subset T^*M$ be the fiberwise starshaped hypersurface
corresponding to the cooriented contact manifold $(S^*M,\ga)$,
and (up to the time change $t \mapsto 2t)$
view the Reeb flow $\gf_\ga^t$ on $(S^*M,\ga)$ as the restriction of the
Hamiltonian flow $\gf_H^t$ on~$T^*M$, where $\Sigma = H^{-1}(1)$ and $H$ is fiberwise homogeneous of degree~$2$
(and smoothened to zero near the zero-section).
For $q \in M$ consider $\Sigma_q = \Sigma \cap T_q^*M$
and the bounded component~$D_q$ of $T_q^*M \setminus \Sigma_q$.
The ``spheres'' $\Sigma_q$ are Legendrian and the ``discs'' $D_q$ are Lagrangian.
Now fix $p \in M$ and let $q$ be such that $p,q$ are non-conjugate.
``Sandwich'' the Hamiltonian~$H$ between a smaller and a larger Riemannian Hamiltonian~$G$ and~$c \2 G$
for a suitable constant $c>0$.
Then continuation maps in Lagrangian Floer homology, the Abbondandolo--Schwarz isomorphism
from the action-filtered Lagrangian Floer homology $HF^n(\gf_G^n(T_p^*M),T_q^*M;\FF)$ of the pair 
of Lagrangian submanifolds $\gf_G^n(T_p^*M), T_q^*M$ 
to the energy-filtered homology $H^{n^2}(\Omega_{pq}M;\FF)$,
and Theorem~\ref{t:main} imply that $\dim HF^n(\gf_H^n(T_p^*M),T_q^*M;\FF)$
grows at least like $e^{\sqrt n}$, uniformly in~$q$.
Since the chain complex of Lagrangian Floer homology in generated by the intersection points 
of the two Lagrangians, 
one concludes that the number of intersections of $\gf_H^n (D_p)$ and~$D_q$
grows at least like~$e^{\sqrt n}$, uniformly in~$q$ for almost every $q \in M$.
Hence the volume of $\gf_H^n (D_p)$ grows at least like $e^{\sqrt n}$.
We refer to Section~4 of~\cite{MacSch11} for details.
Finally, the volume of $\gf_H^n (\Sigma_p)$ also grows like $e^{\sqrt n}$ in view of 
Proposition~4.3 in~\cite{FraLabSch13}.
\proofend


\end{document}